%
\magnification=\magstep1


\def\hexnumber#1{\ifcase#1 0\or1\or2\or3\or4\or5\or6\or7\or8\or9\or
	A\or B\or C\or D\or E\or F\fi }

\font\teneuf=eufm10
\font\seveneuf=eufm7
\font\fiveeuf=eufm5
\newfam\euffam
\textfont\euffam=\teneuf
\scriptfont\euffam=\seveneuf
\scriptscriptfont\euffam=\fiveeuf


\font\tenmsx=msam10
\font\sevenmsx=msam7
\font\fivemsx=msam5
\font\tenmsy=msbm10
\font\sevenmsy=msbm7
\font\fivemsy=msbm5
\newfam\msxfam
\newfam\msyfam
\textfont\msxfam=\tenmsx  \scriptfont\msxfam=\sevenmsx
  \scriptscriptfont\msxfam=\fivemsx
\textfont\msyfam=\tenmsy  \scriptfont\msyfam=\sevenmsy
  \scriptscriptfont\msyfam=\fivemsy
\edef\msx{\hexnumber\msxfam}

\mathchardef\upharpoonright="0\msx16

\def\qed{{\vcenter{\hrule height.4pt \hbox{\vrule width.4pt height5pt
 \kern5pt \vrule width.4pt} \hrule height.4pt}}}
\def\notin{{\in}\kern-5.5pt / \kern1pt}
\def\ok{\vbox{\hrule height 8pt width 8pt depth -7.4pt
    \hbox{\vrule width 0.6pt height 7.4pt \kern 7.4pt \vrule width 0.6pt height 7.4pt}
    \hrule height 0.6pt width 8pt}}
\def\nt{{\leq}\kern-1.5pt \vrule height 6.5pt width.8pt depth-0.5pt \kern 1pt}
\def\sd{{\times}\kern-2pt \vrule height 5pt width.6pt depth0pt \kern1pt}
\def\zp#1{{\hochss Y}\kern-3pt$_{#1}$\kern-1pt}

\def\la{{\langle}}

\def \o {\omega }

\overfullrule=0pt
\openup1.5\jot

\def\square{\qed }
\def \al {\alpha } \def \la {\lambda }
\def \l {\langle } \def \r {\rangle } \def \c {{\bf c}} \def \b {{\bf b}}
\font\titel=cmbx10 scaled \magstep1
\font\utitel=cmbx10 scaled \magstephalf
\font\htitel=cmbx10 scaled \magstep2
\font\namen=cmr10 scaled \magstep2
\font\and=cmr10 scaled \magstep1

{\centerline {\htitel On Gross Spaces}}\footnote{}{{\it 1991
}}

\bigskip

\vskip 2.5 true cm

\centerline{{\namen Saharon Shelah }\footnote{$^1$}{The authors are
supported by the Basic Research Foundation of the Israel Academy of
Science} \footnote{$^2$}{Publication number 468} {\and and } {\namen
Otmar Spinas }$^1$}

\vskip 2.5 true cm

{{\narrower ABSTRACT:  A {\it Gross space} is a vector space $E$ of
infinite dimension over
some field $F$, which is endowed with a symmetric bilinear form $\Phi
:E^2 \rightarrow F$ and has the property that every infinite
dimensional subspace $U\subseteq E$ satisfies dim$U^\perp < $ dim$E$.
Gross spaces over uncountable fields exist (in certain dimensions)
(see [G/O]).
The existence of a Gross space over countable or finite fields (in a
fixed dimension not above the continuum) is independent of the axioms
of ZFC. This was shown in [B/G], [B/Sp] and [Sp2]. Here we continue the
investigation of Gross spaces. Among other things we show that if the
cardinal invariant {\bf b} equals $\omega _1$ a Gross space in
dimension $\omega _1$ exists over every infinite field, and that it is
consistent that Gross spaces exist over every infinite field but not
over any finite field. We also generalize the notion of a Gross space
and construct generalized Gross spaces in ZFC.

}}

\bigskip \bigskip

{\titel 0  Introduction}

\bigskip \bigskip

Let $E$ be a vector space of infinite dimension over some field $F$,
and let $\Phi :E\times E\rightarrow F$ be a symmetric bilinear form.
If $U$ is a subspace of $E$, the {\it orthogonal} complement is denoted by
$U^\perp $. If $U$ has finite dimension, then clearly we have
codim$_EU^\perp \leq $ dim$U$, where codim$_EU^\perp $ is the
dimension of some {\it linear} complement of $U^\perp $ in $E$. As is
well-known from Hilbert space, this is completely false if $U$ is of
infinite dimension. 

In [G/O], the investigation of quadratic spaces sharing the following
strong property has been started:

\bigskip

$(\ast )$  for all subspaces $U\subseteq E$ of infinite dimension:
dim$U^\perp <$ dim$E$.

\bigskip

Such a space we call a {\it Gross space}. In [G/O], the motivation
for this was that Gross spaces are natural candidates to have a small
orthogonal group in the sense that every isometry is the product of
finitely many hyperplane reflections. If a Gross space is the
orthogonal sum of two subspaces, then one of them must be finite
dimensional. Hence, a Gross space is far from having an orthogonal
basis, and so its dimension must be uncountable.

In [G/O], over every uncountable field $F$ a {\it strong Gross space}
has been constructed, i.e. a space sharing the following stronger
property:

\bigskip

$(\ast \ast )$  for all subspaces $U\subseteq E$ of infinite dimension:
dim$U^\perp \leq \aleph _0$.

\bigskip

\noindent Such a construction has been achieved in every uncountable
dimension less or equal the cardinality of the field.   

To construct a Gross space gets the more difficult the smaller the
cardinality of the field is compared with the dimension of the space.

In [B/G], a Gross space of dimension $\aleph _1$ (so Gross$=$strongly
Gross) has been constructed over every countable or finite field. But
for this the Continuum Hypothesis (CH) has been assumed, and hence the
question was raised whether CH is necessary or a construction in ZFC
is possible.

In [Sp1], [B/Sp] and [Sp2], it turned out that this question leads
into set theory, and that in fact it is independent of ZFC.

Here we answer several open questions from these papers, and we also
continue the investigation of Gross spaces over uncountable fields in
[G/O]. 



A difficult result in [B/Sp] says that if ${\bf b}=\omega _1$, then a Gross
space of dimension $\aleph _1$ exists over every field which is the
extension of some finite or countable field by countably many
transcendentals. Here ${\bf b}$ is defined as the
minimal cardinality of a family of functions from the natural
numbers to themselves which is unbounded under eventual dominance.
Clearly $\omega _1\leq {\bf b}\leq {\bf c}$. But ZFC does not decide
where ${\bf b}$ lies exactly. In
$\S 2$ we will show that if ${\bf b}=\omega _1$, then a Gross space of
dimension $\aleph _1$ exists over {\it every} infinite field.  

In $\S 3$ we show that ${\bf b}=\omega _1$ may hold but no symmetric
bilinear space of dimension $\aleph _1$ over any finite field is a
Gross space. We prove three variants of this. The first uses a forcing
\relax from [Sp2] which kills Gross spaces, the second involves the splitting number ${\bf s}$, and the
third uses a model from [B/Sh] where simultaneously $P_{\aleph _1}$-
and $P_{\aleph _2}$-points exist. In this model there exist no Gross
spaces over finite fields at all. 

For a couple of years the main open problem about Gross spaces has been
whether there exists a ZFC
model where there exists no Gross space over any countable or finite
field in any dimension. The first author thinks that he has constructed a
model for this, using a new iteration technique. The paper has not yet 
been written.

In $\S 4$ we investigate a natural generalization of the Gross
property:

\bigskip  

$(\ast \ast \ast )$  for all subspaces $U\subseteq E$ of dimension
$\geq \lambda $:
dim$U^\perp <$ dim$E$

\bigskip

\noindent Here $\lambda $ is an infinite cardinal less or equal the dimension
of $E$. A
space sharing the property $(\ast \ast \ast )$ is called a $\lambda
$-{\it Gross space}. Hence a Gross space is an $\omega $-Gross space. 

We concentrate on results in ZFC. We show
that $\omega _1$-Gross spaces of dimension $\aleph _1$ and $\aleph _2$
always exist over every countably infinite field. We also show that
$|F|^\lambda $ is an upper bound for the dimension of a $\lambda
$-Gross space over $F$. Hence, in ZFC these results are maximal. We do not know
whether these results hold also for finite fields. But we show that
if $\la >\o _1 $ is
regular uncountable, then a $\la ^+$-Gross space of dimension  $\la
^+$ exists over any field.

We also investigate the situation for uncountable fields. The
construction in [G/O] yields strong Gross spaces of dimension at most
as large as the cardinality of the base field. Is it possible to
enlarge the dimension of the space while keeping the size of the field
small? Over every field of uncountable cardinality $\lambda $ we
construct a $\lambda $-Gross space of dimension $\lambda ^+$, as well
as a $\lambda ^+$-Gross space of dimension $\lambda ^{++}$. Again, by
the upper bound for the dimension of a $\lambda -$Gross space
mentioned above, these results are maximal in ZFC.

We remark that all the investigations mentioned above may take place
in a rather more general context. First, everything remains true if we
skip to orthosymmetric sesquilinear forms over a division ring endowed
with an involutory antiautomorphism (see [G] for the definitions).
Second, using a representation theorem
for AC-lattices equipped with a polarity,
the results can be transferred to the level of abstract
ortho-lattices. This has been announced in [Sp3].

For a survey on the whole subject refer to [Sp4].

\bigskip \bigskip

{\titel 1 Notation and definitions}

\bigskip \bigskip

{\utitel 1.0 Forms}

\bigskip

Let $F$ be a commutative field of arbitrary characteristic. Let $E$ be
a vector space over $F$, endowed with a symmetric bilinear form
$\Phi :E\times E\rightarrow F$, i.e. $\Phi $ is linear in both
arguments and $\Phi (x,y)=\Phi (y,x)$ always. Most
of the spaces we construct
will have {\it isotropic} vectors, i.e. nonzero vectors $x$ such that
$\Phi (x,x)=0$. For a subspace $U\subseteq E$ the {\it orthogonal}
complement $U^\perp $ is the subspace
$\{ x\in E:\forall y \in U$ $\Phi (x,y)=0\} $. We call $\l E,\Phi \r $
{\it nondegenerate} if $E^\perp =\{ 0\} $.
By $E^*$ we denote the $F-$vectorspace of linear functionals
$E\rightarrow F$. If in $E$, we have fixed a basis $\langle e_\alpha :
\alpha \in
I\rangle $, then for a vector $x\in E$ the {\it support} of $x$,
denoted by supp$(x)$, is the unique finite set of $\alpha $'s such that, in
the representation of $x$ by $\langle e_\alpha :\alpha \in
I\rangle $, $e_\alpha $ has a nonzero coefficient.

\bigskip  

{\utitel 1.1 Set theory}

\bigskip 

For the theory of forcing refer to [B], [J], [K] or [Sh1].

For sets $A, B$, the set of functions from $A$ to $B$ is denoted by
${{^A}B}$. If $\kappa ,\lambda $ are cardinals, then $\kappa ^\lambda $
denotes the cardinality of the set ${ ^\lambda }\kappa $, and $[A]^\lambda $
denotes the
set of all subsets of $A$ which have cardinality $\lambda $. By ${\bf
c}$ we denote the cardinality of the continuum. If $A$ is
well-ordered, then o.t.$(A)$ is its order type. 

If $f$, $g\in \, {{^\o }\o } $, then we say $f$ {\it eventually dominates}
$g$, and we write $g< ^\ast f$, if $\exists k\forall n\geq k \,
g(n)< f(n)$. A family $F$ of members of ${{^\o } \o }$ is $<^*$-{\it unbounded}
if there is no $f\in \, {{^\o }\o }$ such that $\forall g\in $ ran$F\, g<^\ast
f$. $F$ is called {\it dominating}, if every function from ${^ \o }\o $
is eventually dominated by some member of $F$. 

Then $\b $ is defined as the minimal cardinality of a $<^*
$-unbounded
family in ${^\o } \o $. 

The cardinal invariant ${\bf d}$ is defined as the minimal cardinality
of a dominating family.

For $A$, $B\in [\o ]^\o $, we write $A\subseteq ^\ast B$ if
$A\setminus B$ is finite. 
A family ${\cal S}$ of members of $ [\omega ]^\o $ is called a {\it splitting
family}, if for every $A\in [\o ]^\o $ there exists $S\in $ ran${\cal S}$
such that $A\cap S$ and $A\setminus S$ are both infinite. 

The splitting number
${\bf s}$ is defined as the minimal cardinality of a splitting
family.

Here ${\bf b}$, ${\bf d}$ and ${\bf s}$ are instances of so called cardinal
invariants of the continuum which have found lots of applications
throughout mathematics. See [vD] for an introduction.

The following theorems form part of the folklore in set theory. Proofs of the first two of them may be found in
[vD].

\bigskip

{\bf Theorem 1}. {\it ${\bf b}$, ${\bf d}$ and ${\bf s}$ are uncountable
cardinals $\leq {\bf c}$; ${\bf b}, {\bf s}\leq {\bf d}$, and $\b $ is regular.}

\bigskip

The axioms of ZFC do not decide where exactly ${\bf b}$, ${\bf d}$ and
${\bf s}$
lie.

\bigskip

{\bf Theorem 2}.  {\it Let $\kappa $ and $\lambda $
be regular cardinals with $\o _1\leq \kappa \leq \lambda $. There
exists a forcing preserving all cardinals such that in the extension
obtained by forcing with it, ${\bf c} =\lambda $ and ${\bf b}= {\bf
d}={\bf s} = \kappa $ hold.}

\bigskip

{\bf Theorem 3 [Sh2]/[BS].}  {\it $\omega _1 ={\bf b} <{\bf s}$
is consistent with} ZFC {\it as well as is}
{\it $\omega _1 ={\bf s}<{\bf b} $.}





\def \la {\lambda }  

\bigskip \bigskip

{\titel 2  When ${\bf b}=\o _1 $, a Gross space exists over every infinite
field}

\bigskip \bigskip

In [B/Sp, $\S 4$] it has been shown that the assumption ${\bf b}=\omega _1$
implies that a Gross space of dimension $\o _1$ exists over any field
which is the extension of an arbitrary
finite or countable field by countably many transcendentals. 
Here we show that ${\bf b}=\o _1$
is enough to construct a Gross space of dimension $\o _1$ over any
countably infinite field. In $\S 3$, we will show that ${\bf b}=\o _1$ may hold
but no Gross spaces exist over finite fields.

The construction will use the filtrations given by the following
Lemma.

\bigskip

{\bf Lemma.}  {\it There exists a family $\langle A^\alpha _k :
\alpha <\omega _1, k<\omega \rangle $ of finite sets such that
for any $\al, \beta <\o _1 $ and $k,l <\o $ the following requirements
are satisfied:

(1) $\al =\bigcup_{k<\o }A^\al _k$

(2) if $k<l$ then $A^\al _k \subseteq A^\al _l$

(3) if $\beta <\al $ and $\beta \in A^\al _k$ then $A^\beta _k
=A^\al _k\cap \beta $}

\bigskip

{\bf Proof.}  We define $A_k^\alpha $ by
induction on $\alpha $. 

For $\alpha =0$ let $A_k^\alpha =\emptyset $.

For $\alpha =\beta +1$ let $A_k^\alpha =A_k^\beta \cup \{ \beta \} $

For $\alpha $ a limit let $\langle \alpha _n :n<
\omega \rangle $ be an increasing sequence such that $\alpha
=\sup_{n<\omega }\alpha _n$. Now choose $\langle k_n : n<\omega
\rangle $ increasing such that for all $n$ we have $\alpha _n\in
A_{k_n}^{\alpha _{n+1}}$. Define 

$$A_m^\alpha =A_m^{\alpha _0}$$

\noindent for $m\leq k_0$ and

$$A_m^\alpha =A_m^{\alpha _{l+1}}$$

\noindent for $k_{l}< m\leq k_{l+1}$. It is easy to check that this
works. \hfill $\square $ 

\bigskip

{\bf Theorem.}  {\it Assume} ${\bf b}=\omega _1$. {\it Then a
Gross space of dimension $\aleph _1$ exists over every infinite field.}

\bigskip

{\bf Proof.}  Let $F$ be a countable, infinite field, and let $E$ be a
vector space over $F$, spanned by a basis $\langle e_\alpha : \alpha <
\omega _1\rangle $. 

Let $\langle p_n :n\in \omega \rangle $ enumerate all polynomials in
finitely many variables with coefficients in $F$.

In $F$ choose $\langle a_n: n<\omega \rangle $
such that for any polynomial $p[X_1,\dots ,X_k]$, for any $n$ larger than some integer depending on $p$
we have $a_n \not\in \{ p(a_{l_1},\dots ,a_{l_k}): l_1,\dots ,l_k <n\}
$. 

Finally, let $\langle f_\alpha :\alpha <\omega
_1\rangle $ be an unbounded, well-ordered family of increasing functions 
in $\omega ^\omega $.

For any $\al <\omega _1$ we define a function $h_\al :\al \rightarrow
\omega $ by induction, using the filtrations from the Lemma. Assume
$\gamma \in A_{k+1}^\al\setminus A_k^\al $ and assume $h_\beta $ for $\beta
<\al $, $h_\al \upharpoonright A^\al _k $ and $h_\al
\upharpoonright (A^\al _{k+1} \cap \gamma )$ have been defined.

Let $h_\al (\gamma )=n $ such that
$a_n$ is distinct from all the
finitely many scalars of the form 
$$p_i(a_{l_1}, a_{l_2},\dots )$$
where $i< f_\al (k+1)$ and any $l_j$ belongs to the range of one of
the following functions: $h_\al \upharpoonright A^\al _k$, $h_\al
\upharpoonright (A^\al _{k+1}\cap \gamma )$, $h_\beta \upharpoonright
A^\beta _{k+1}$ for $\beta \in A^\al _{k+1}$.

Now define a symmetric bilinear form 
$\Phi :E\times E \rightarrow F$ as follows: For $\al <\beta  <\omega
_1$ define:

$$\Phi (e_\al , e_\beta )=a_n \hbox{ if and only if }h_\beta (\al ) =n
$$

$ \Phi (e_\beta, e_\beta ) $ may be defined arbitrarily. 

We have to show that $\l E, \Phi \r $ is Gross. Assume that this is not
true. So there is a subspace $U$ of infinite dimension, spanned by
a basis $\langle y_k : k <\o \rangle $, such that dim$U^\perp $ is
uncountable. Each $y_k$ has a representation

$$y_k =\sum_{l=1}^{m_k}b_{kl} e_{\al (k,l)}$$

\noindent with nonzero coefficients $b_{kl}$. 
Choose $\al ^* <\omega _1$ such that $U\subseteq $
span$\langle e_\al ,\al <\al ^*\rangle $. Using a $\Delta -$system
argument and the bilinearity of the form, in $U^\perp $ we may certainly find 
vectors $\langle z_\iota ,\iota <
\omega _1\rangle $, all of them having the same nonzero
coefficients
in their representation, say

$$z_\iota =\sum_{l=1}^n c_l e_{\beta (\iota ,l)}$$ 

\noindent
such that for all $\iota _1<\iota _2<\omega _1$ we have $\al ^*< \beta (\iota
_1 ,1)<\dots <\beta (\iota _1, n)<\beta (\iota _2, 1)<\dots <\beta
(\iota _2,n)$. Furthermore, we may assume that for some $k^*<\omega $,
for all $\iota <\omega _1$, $1\leq l_1<l_2\leq n$ and $1\leq l\leq
n$ we have

$$\beta (\iota ,l_1)\in A^{\beta (\iota ,l_2)}_{k^*} \hbox{ and }
\alpha ^*\in A^{\beta (\iota ,l)}_{k^*}\eqno(1)$$

Since $\langle f_\alpha ,\al < \omega _1 \rangle $ is unbounded and
well-ordered there
exists $j<\o $ such that $\langle f_{\beta (\iota ,n)}(j), \iota \in \o _1 \rangle $ is
unbounded in $\o $. Since the $f_\al $ are increasing this remains true if $j$ is replaced by
any $j'\geq j$. 

Now choose $k<\o $ such that the minimal $j'$ with the property

$$\{ \al (k,1), \dots ,\al (k,m_k)\} \subseteq A^{\al ^*}_{j'+1}\eqno(2)$$

\noindent is larger than max$\{ j, k^*\} $. Note that by (1) and
property (3) of the filtrations from the Lemma, $j'$ is 
minimal such that (2) holds for any $\beta (\iota
,l)$ instead of $\al ^*$. 

Choose $l_0\in \{ 1, \dots , m_k\} $ such that $\al (k,l_0)$ is maximal in
$A^{\al ^*}_{j'+1} \setminus A^{\al ^*}_{j'}$. Then the polynomial

$$-\sum_{\langle l,l'\rangle \in \{1,\dots , n\} \times \{ 1,\dots ,
m_k\} \setminus \langle n,l_0 \rangle } 
b_{k l_0}^{-1}c_n^{-1}b_{kl'}c_lX_{ll'}$$

\noindent has its number, say $i$, in the enumeration of all
polynomials fixed in the beginning. There exists $\iota <\omega _1$
such that $f_{\beta (\iota ,n)}(j'+1)>i$. From the definition of $\Phi
$, using (1), we conclude 

$$\Phi (e_{\beta (\iota ,n)}, e_{\al (k,l_0)})\ne -\sum_{\langle
l,l'\rangle \in \{
1,\dots ,n \} \times \{ 1, \dots , m_k\} \setminus \langle
n,l_0\rangle }
b_{kl_0}^{-1}c_n^{-1}
b_{kl'}c_l\Phi (e_{\beta (\iota ,l)}, e_{\al (k,l')})$$

\noindent But then $\Phi (z_\iota ,y_k)\ne 0$, a contradiction.\hfill
$\square $

\bigskip \bigskip

{\titel 3  Gross spaces over finite fields}

\bigskip  \bigskip

The first Lemma shows that $\c $ is an upper bound for the dimension
of a Gross space over finite or countable fields.

\bigskip

{\bf Lemma 1.}  {\it Let $\l E,\Phi \r $ be a symmetric bilinear space
over $F$, and let $\la $ be an infinite cardinal $\leq $ the dimension
of $E$. If 
$\l E,\Phi \r $ has the property}

$$(\ast \ast \ast ) \hbox{ for all subspaces }U\subseteq E\hbox{ of
dimension }\geq \la \hbox{: dim}U^\perp <\hbox{ dim}E$$

\noindent {\it then} dim$E\leq |F|^\lambda $. 

\bigskip

{\bf Proof.}  Let dim$E=\kappa $ and let $\langle e_\alpha :\alpha <
\kappa \rangle $ be a basis. For $\alpha <\kappa $ define

$$f_\alpha : \lambda \rightarrow F, \, \beta \mapsto \Phi (e_\alpha
,e_\beta )$$

We get $\kappa =\bigcup_{f\in ^\lambda F}A_f$, where $A_f=\{ \alpha
<\kappa : f_\alpha =f\} $ and hence $\kappa =\sum_{f \in ^\lambda
F}|A_f|$. 

So if $\kappa >|F|^\lambda $ there exists $B\subseteq
\, {^\lambda F}$ such that $\sum_{f\in B}|A_f| =\kappa $ and $|A_f|\geq 2$,
for all $f\in B$. Fix $\alpha _f \in A_f$ for $f\in B$. We conclude
that 

$$\bigcup_{f\in B} \{ e_\alpha -e_{\alpha _f}: \alpha \in A_f \hbox{
and } \alpha \ne \alpha _f \} $$

\noindent is linearly independent, of cardinality $\kappa $ and a
subset of $\langle e_\beta : \beta <\lambda \rangle ^\perp $. \hfill
$\square $

\bigskip

In [Sp2], the following forcing to kill Gross spaces was introduced.

Let $F$ be a field. The direct sum $\oplus_{n<\o }F $ of $\o $ copies
of $F$ naturally bears a vector space structure over $F$

The forcing $P^F$ consists of pairs $\l
s,A\r $ where 

(1) $s=\l s_0 ,\ldots ,s_k \r $ is a finite sequence of vectors
$s_i\in \oplus_{n<\o } F$

(2) $A\subseteq (\oplus_{n<\o } F)^*$ is a finite set of linear functionals.

The ordering on $P^F$ is defined as follows:

$$\halign{ #&#\hfill \cr
$\l s',A'\r \hbox{ extends } \l s,A\r $ if and only if &$s'\supseteq s$ and
$A'\supseteq A$ and \cr &$\forall s_i'\in $ ran$(s')$ $\setminus $
ran$(s)$ $\forall f\in A$ $f(s_i')=0$\cr }$$

It is easily seen that two conditions with the same first coordinate
are compatibel. Hence, if $F$ is countable or finite $P^F$ has the
countable chain conditon. 

Let $\l E,\Phi \r $ be an inner product space over $F$ and $U$ an
arbitrary subspace of infinite dimension. Then forcing with $P^F$
introduces in $U$ an infinite set of linearly independent vectors such
that every vector in $E$ is orthogonal on all but finitely many of
them. Moreover, iterating forcing with $P^F$ such that $F$ again and again runs
through all finite and countable fields produces a model where Gross
spaces of any dimension less than the length of the iteration (which
is supposed to be a regular cardinal) do not
exist (see [Sp2] for the proofs). 

By the Theorem in $\S 2$ above, if ${\bf b}=\o _1$ holds a Gross space of
dimension $\o
_1$ exists over every infinite field. Hence if $F$ is infinite, the
forcing $P^F$ which kills Gross spaces over $F$ adds a dominating
function (see [Sp2]). Here we show that in case $F$ is finite $P^F$
does not destroy unbounded families, more exactly it is almost ${^\o
}\o -$bounding. Hence, by a preservation theorem
in [Sh1] or [Sh2], iterating forcings $P^F$
where $F$ runs again and again through all finite fields yields a
model where a Gross space in dimension $\o _1$ exists over every
infinite field but not over any finite field, provided that we start
with a ground model satisfying ${\bf b}=\o _1$.

\bigskip

{\bf Lemma 2.}  {\it Assume that $F$ is a finite field and $\tau $ is a
$P^F$-name for a function in ${^\o }\o $. Then there exists $f\in {^\o }\o
$ in the ground model such that for every $p\in P^F$ there are only
finitely many
$n$ such that $p\vdash \tau (n)>f(n)$.}

\bigskip

{\bf Proof.}  At first, let us introduce the following notation. If
$s=\langle s_0, \ldots ,s_k\rangle $ is a sequence of vectors in
$\oplus_{n<\o }F$ and $\phi $ is a formula in the language of forcing with
$P^F$ we write $s\vdash ^* \phi $ if and only if there exists a finite
$A\subset (\oplus_{n<\o }F)^* $ such that $\langle s,A\rangle \vdash \phi $.

By induction on $i<\o $ we will define an increasing sequence $\langle
n_i: i<\o \rangle $ such that:

(1) $n_0=0$ and

(2) for every $i$, if $s=\langle s_0,\ldots ,s_k\rangle $ is a
sequence of vectors in $\oplus _{j<n_i}F$ and $f_0,\ldots ,f_{n_i}$ are linear
functionals on $\oplus_{j<n_{i+1}}F$, then there is a sequence $s'$ of
vectors in $\oplus_{j<n_{i+1}}F$ which extends $s$, such that ran$(s')$
$\setminus $ ran$(s)\subseteq \bigcap_{k\leq n_i}$ker$f_k$ and for
some $\sigma \in {^{n_i}\o }$ we have 
$s'\vdash ^* \tau \upharpoonright n_i =\sigma $.

We have to show that for given $n_i$ there exists $n_{i+1}>n_i$ such
that (2) holds. Assume that this is not true. So we may find $s$ with
ran$(s) \subseteq \oplus_{j<n_i}F$ such that for infinitely many $j>n_i$
there are $f^j_0,\ldots ,f^j_{n_i} \in (\oplus_{u<j}F)^*$ such that there is
no $s'$ as in the conclusion of (2). 

Now every $(\oplus_{u<j}F)^*$ is finite.
Hence by K\"onig's Lemma, there exist $f_0,\ldots ,f_{n_i}\in (\oplus_{j<\o
}F)^*$ such that for some $B\in [\o ]^\o $ we have $\forall j \in B\exists j'>j\forall k\leq n_i$
$f_k\upharpoonright \oplus_{u<j}F=f_k^{j'}\upharpoonright \oplus_{u<j}F$. 

Since $\langle s,\{ f_0,\ldots ,f_{n_i}\} \rangle $ is a condition in
$P^F$ there exists a stronger condition $\langle s',A'\rangle $ which
decides the first $n_i$ values of $\tau $. Choose $j\in B$ so large
that ran$(s')\subseteq \oplus_{u<j}F$. Now for some $j'>\max \{ j,n_i\} $
we know $\forall k\leq n_i$ $f_k\upharpoonright \oplus_{u<j}F =
f_l^{j'}\upharpoonright \oplus_{u<j}F$. But then $s'$ satisfies the
conclusion of (2) for $s$ and $f_0^{j'},\ldots
,f_{n_i}^{j'} $, a contradiction. 

Hence we have shown that $\langle
n_i:i<\o \rangle $ can be chosen increasing and satisfying (1) and (2).

Next we define $f:\o \rightarrow \o $ as follows: If $n_i\leq
n<n_{i+1}$, then 

$$f(n)=\max \{ \sigma (n): \hbox{ there exists }s \hbox{ such that ran}
(s)\subseteq
\oplus_{j<n_{i+2}}F \hbox{ and }s\vdash ^* \tau \upharpoonright
n_{i+1}= \sigma \} $$

We claim that $f$ is as desired. So let $\langle s,\{ f_0,\ldots
,f_m\} \rangle \in P^F$ and let $i<\o $ such that $n_i>m$ and
ran$(s)\subseteq \oplus_{j<n_i}F$. Assume now $n\geq n_i$. We will find an
extension of $\langle s,\{ f_0,\ldots ,f_m \} \rangle $ forcing
``$f(n)\geq \tau (n)$''. 

Choose $j$ such that $n_j\leq n<n_{j+1}$. By
construction, there is $s'$ extending $s$ such that
ran$(s')\subseteq \oplus_{u<n_{j+2}}F$, ran$(s')$ $\setminus $
ran$(s)\subseteq \bigcap_{k=1}^m $ker$f_k$ and $s'\vdash ^* \tau
\upharpoonright n_{j+1}=\sigma  $, for some $\sigma $. Choose a finite
$A\subseteq (\oplus_{u<\o }F )^*$ containing $\{ f_0,\ldots ,f_m\} $ such that
$\langle s',A\rangle \vdash \tau \upharpoonright n_{j+1}
=\sigma$. Then $\langle s', A\rangle $ extends $\langle s, \{
f_0,\ldots ,f_m\} \rangle $ and $\langle s', A\rangle \vdash
f(n)\geq \tau (n)$. \hfill $\square $

\bigskip

{\bf Remark:}  It is easy to see that Lemma 2 says much more than that
$P^F$ is almost ${^\o }\o -$bounding (this notion is defined in
[Sh2]). In [Sh2] it is proved that this property is preserved under
countable support iterations. In the new edition of [Sh1] it is proved
that this is true also for finite support iterations. A more accessible proof
of this preservation theorem (for the finite support case) can be found in
[Go, Example 8.5]. It is not
difficult to see that $P^F$ adds a Cohen real, hence it is not ${^\o
}\o -$bounding.

\bigskip

{\bf Corollary.}  {\it Assume $P$ is a finite or countable support
iteration of forcings $P^F$ where $F$ is a finite field. If in $V$,
$\langle f_\alpha :\alpha <\kappa \rangle $ is an unbounded family of
functions in $\o ^\o $, then $\langle f_\alpha :\alpha <\kappa \rangle
$ is unbounded in $V^P$.}

\bigskip

{\bf Proof.}  Preservation of unboundedness at successor steps of the
iteration follows from Lemma 1: Assume $\tau $ is a $P$-name for a
function in $\o ^\o $ such that some $p\in P$ forces ``$\tau \hbox{
bounds } \langle f_\alpha :\alpha <\kappa \rangle $''. Choose $f\in \o
^\o \cap V$ for
$\tau $ as in the Lemma. Since $\langle f_\alpha :\alpha <\kappa
\rangle $ is unbounded in $V$ there exists $\alpha $ such that for
infinitely many $n$ we have $f(n)<f_\alpha (n)$. Find $p'\in P$
extending $p$ and $n_0$ such that $p'\vdash \forall n\geq n_0 \,
\tau (n)>f_{\al }(n)$. But then for infinitely many $n$ we would have
$p'\vdash \tau (n)>f(n)$, a contradiction. 

The limit steps are handled by a general preservation theorem due to
the first author (see [Sh2] or [Sh1-new ed.,VI$\S 3$]).\hfill $\square
$ 

\bigskip

\relax From the Corollary, the Theorem in $\S 2$ and the results about
forcing with
$P^F$'s in [Sp2] we obtain the following theorem:

\bigskip

{\bf Theorem 1.}  {\it Assume that the ground model $V$ satisfies} ${\bf
b}=\o _1$. {\it Let $P= \lim \langle P_\alpha ,Q_\alpha :\alpha
<\kappa \rangle $ be a finite or countable support iteration of length
$\kappa $ where} $\kappa >\o _1$ {\it is regular such that 

(1) each $Q_\alpha $ is $P^F$ defined in $V^{P_\alpha }$ where $F$ is
a finite field, and

(2) for each finite field $F$, cofinally many times we have $Q_\alpha
=P^F$.

\noindent Then in $V^P$, there exists a Gross space of dimension $\o
_1$ over every infinite field but there is no Gross space of dimension
$\kappa $ over any finite field.}

\bigskip

Another way to get a model for the conclusion of the Theorem above is
shown by the following Lemma combined with Theorem 3, $\S 1$, and the
Theorem in $\S 2$.

\bigskip

{\bf Lemma 3.}  {\it Assume that $\l E,\Phi \r $ is a symmetric
bilinear space over a finite field and the dimension of} $E$ {\it has
uncountable cofinality. If} ${\bf s}>$ dim$E$ {\it then $\langle E,\Phi \rangle $ is not
Gross.}  

\bigskip

{\bf Proof.}  Let dim$E=\kappa $. In $E$, choose a basis 
$\langle e_\alpha ,\alpha <\kappa
\rangle $ and let $\langle y_n,n<\omega \rangle $ span a subspace. For
every $\alpha <\kappa $ let 

$$f_\alpha :\omega \rightarrow F, \, n \mapsto \Phi (e_\alpha ,y_n)$$

\noindent
where $F$ is the base field. Let ${\cal A}= \langle f_\alpha ^{-1}\{
a\} :\alpha <\kappa , a\in F\rangle $. Since ${\bf s}>\kappa $, ${\cal
A}$ is not a splitting family. Hence there is an infinite $A\subseteq
\omega $ such that for all $B\in $ ran${\cal A}$, either $A\subseteq ^*B$ or
$A\cap B$ is finite. For any $\alpha <\kappa $, $\langle f_\alpha
^{-1}\{ a\} , a \in F\rangle $ is a finite partition of $\omega $.
So there exist $a_\alpha $ and $n_\alpha $ such that $A\setminus
n_\alpha \subseteq f_\alpha ^{-1}\{ a_\alpha \} $. By cf$(\kappa
)>\omega $ there exist $a$, $n$ and $X\in [\kappa ]^\kappa $ such
that for every $\alpha \in X$ we have $a_\alpha =a$ and $n_\alpha =n$.
Now it is easy to see that every vector $e_\alpha -e_\beta $ with
$\alpha $, $\beta  \in X$ is in the orthogonal complement of
span$\langle y_m , m\in A\setminus n\rangle $.\hfill $\square $ 

\bigskip

The model from [B/Sh] where simultaneously simple $P_{\aleph _1}$- and
$P_{\aleph _2}$-points exist shows that a stronger statement than the
conclusion of Theorem 1 is true. A $P_\kappa $-point, where $\kappa $
is a regular uncountable cardinal, is a filter on $\o $ which is
generated by a tower (i.e. a $\subseteq ^*$-decreasing family in $[\o
]^\o $) of length $\kappa $.

\bigskip

{\bf Theorem 2.}  {\it It is consistent that a Gross space of dimension
$\omega _1 $ exists over every infinite field but no Gross space
exists over any finite field in any dimension.}

\bigskip

{\bf Proof.}  In the model of [B/Sh] where simple $P_{\aleph _1}$- and
$P_{\aleph _2}$-points exist, $\b =\o _1$ holds. Hence by the Theorem
in $\S 2$, a
Gross space of dimension $\aleph _1$ exists over any infinite field. Using
a similar argument as in the proof of Lemma 3, one shows that a
$P_{\aleph _1}$-, $P_{\aleph _2}$-point rules out the existence of a
Gross space over any finite field in dimension $\aleph _1$, $\aleph _2$
respectively. But in that model $\c =\o _2$. Hence by Lemma 1 
we are done.\hfill $\square $

\bigskip

{\bf Question:} Does there exist a ZFC--model where there exists no
Gross space over any finite or countable field in any dimension?

\bigskip

In view of the Theorem in $\S 2$, a natural question to ask is whether
the assumption ${\bf s}=\omega
_1$ is strong enough to construct a Gross space of dimension $\aleph _1$ over
a finite field. The answer is ``no''.

\bigskip

{\bf Theorem 3.}  {\it It is consistent that ${\bf s}=\omega _1$ holds
and there exists no Gross space in dimension $\omega _1$ over any
finite field.}

\bigskip

{\bf Proof.}  Since every finite or countable field can be coded by a
real, it is not difficult to see that the forcing notion $P^F$, where
$F$ is finite or countable, is Souslin in the sense of [J/Sh]; i.e.,
the set of conditions can be viewed as an analytic set of reals,
whereas the ordering and the incompatibility relation are analytic
subsets of the plane. In [J/Sh] it is proved that if we start with a
model satisfying CH, then in any finite support extension of Souslin
forcing notions having the countable chain condition, the reals of the
ground model remain a splitting family. Hence in $V^{P}$, where $P$ is
the finite support iteration of Theorem 1 and $V$ satisfies CH, 
${\bf s}=\o _1$ holds and there are no Gross spaces in any dimension
$<{\bf c}$ of uncountable cofinality.  

\bigskip

Is there another cardinal invariant such that its being $\aleph _1$
implies that Gross spaces over finite fields exist? The largest of the
``classical'' cardinal invariants is ${\bf d}$. The model in [Sh3, $\S 2$]
shows that even ${\bf d }=\o _1$ does not suffice for our purpose. In that
model, if $\l E,\Phi \r $ is a quadratic space spanned by a basis $\l
e_\al :\al <\o _1\r $ over some finite field $F$, then the following
holds: 

\bigskip
{\bf Fact.}  {\it If $\l u_n :n<\o \r $ is a sequence of pairwise
disjoint subsets
of $\o $ such that every $u_n$ has size $|F|^n +1$, then we can find
$\l H_n: n<\o \r $, where each $H_n$ is a family of at most $n$
functions $u_n\rightarrow F$, such that 

$$\forall \al <\o _1 \exists n_\al <\o \forall n\geq n_\al \exists h
\in H_n \forall m\in u_n \,  [\Phi (e_\al ,e_m)=h(m)]$$}

It is not difficult to see that this implies that $\l E,\Phi \r $ is
not Gross (solve many systems of homogeneous equations). 

Next we will show that if there is a family of $\aleph _1$ many
meagre sets of reals which is cofinal with respect to inclusion
in the set of all meagre sets, then a Gross space of dimension
$\aleph _1$ exists over every finite field. We work in $2^\o $
considered as topological product of the discrete space $2=\{ 0,1\} $.
Denote by ${\cal M}$ the set of all meagre sets in $2^\o $. Now
cof$({\cal M})$ is the cardinal invariant defined as follows:

$$\hbox{cof}({\cal M})=\min \{ |{\cal F}| : {\cal F} \hbox{ is a
family of meagre sets such that }\forall A\in {\cal M}\exists B\in
\hbox{ ran}{\cal F} \, A\subseteq B \}$$

\def \be {\beta } \def \la {\lambda } \def \i {\iota } \def \ka {\kappa }
\bigskip

{\bf Theorem 4.}  {\it Assume} cof$({\cal M})=\o _1$. {\it A Gross space of
dimension $\aleph _1$ exists over every finite field.}  

\bigskip

{\bf Proof.}  Using the assumption it is standard work to construct a
family $\l r_\al :\al <\o _1\r $ of reals such that, if $A\subseteq
(2^\o )^n $ is a meagre set for some $n$, then $\exists \al <\o _1
\forall \al <\al _1<\ldots <\al _n \, \l r_{\al _1}, \ldots ,r_{\al
_n}\r \not\in A$.

By induction on $\al <\o _1$ choose one-to-one functions $h_\al :\al
\rightarrow \o $ such that $\o \setminus $ ran$(h_\al )$ is infinite
and for every $\beta <\al $ the set $\{ \gamma <\beta :h_\beta (\gamma
) \ne h_\al (\gamma )\} $ is finite. (This is a well-known construction
of an Aronszajn tree due to Todorcevic.)

Now let $F$ be a field and $E$ a vector space over $F$ of dimension
$\aleph _1$, spanned by a basis $\l e_\al :\al <\o _1 \r $. Define a
symmetric bilinear form $\Phi $ on $E$ as follows: For $\al <\beta <\o
_1$ set:

$$\Phi (e_\al ,e_\beta )=r_\beta (h_\beta (\al ))$$

We claim that $\l E,\Phi \r $ is Gross. If this is not true there are
families of vectors $\l y_k :k<\o \r $ and $\l z_\iota :\iota <\o _1\r
$ such that $\Phi (y_k,z_\iota )=0$ always, $\l $dom$(y_k):k<\o \r $
and $\l $dom$(z_\iota ):\iota <\o _1\r $ are families of pairwise
disjoint sets and the sets of the latter have all the same
cardinality and are disjoint from some $\al ^*<\o _1$ such that
span$\l y_k:k<\o _1\r \subseteq $ span$\l e_\al :\al <\al ^*\r $.
Let

$$y_k=\sum_{l=1}^{m_k} b_{kl}e_{\alpha (k,l)}$$

$$z_\iota =\sum_{l=1}^n c_l e_{\beta (\iota, l)}$$

\noindent such that each $m_k>0$ and $n>0$, each $b_{kl}, c_l\ne 0$
and $\be (\i ,1)<\ldots <\be (\i ,n)$ always. 

Define $A\subseteq (2^\o )^n$ as follows:

$$A=\{ \l r_1,\ldots ,r_n\r : \{ k: \sum_{1\leq l\leq m_k}\sum_{1\leq
l'\leq n}b_{kl}c_{l'}r_{l'}(h_{\al ^*}(\al (k,l)))\ne 0\} \hbox{ is
finite}\}$$

It is not difficult to see that $A$ is meagre. Hence by construction we may
choose $\i $ so large that $\l r_{\be (\i ,1)}, \ldots , r_{\be (\i
,n)}\r \not\in A$. Consequently, by the choice of the $h_\al $ we may
find $k$ such that

$$\halign{#& \hfill #\hfill  \cr
$\Phi (y_k,z_\i )=$&$\sum_{1\leq l\leq m_k}\sum_{1\leq l'\leq n}
b_{kl}c_{l'}r_{\be (\i ,l')}(h_{\be (\i ,l')}(\al (k,l)))=$\cr 
&$\sum_{1\leq l\leq m_k}\sum_{1\leq l'\leq n}
b_{kl}c_{l'}r_{\be (\i ,l')}(h_{\al ^*}(\al (k,l))) \ne 0$\cr }$$

This is a contradiction.\hfill $\square $

\bigskip
{\bf Remark.}  Theorem 4 is true if the assumption cof$({\cal M})=\o
_1$ is replaced by cof$({\cal N})=\o _1$ where ${\cal N}$ denotes the
set of all Lebesgue measure zero sets of reals. The same proof works
if we replace ``meagre'' by ``measure zero''. The proof also works for
arbitrary (not necessarily finite) fields. All this is not astonishing
since Cicho\'n's diagram (see [F]) tells us that ${\bf b}\leq {\bf
d}\leq $ cof$({\cal M}) \leq $ cof$({\cal N})$.
  
\bigskip \bigskip

{\titel 4  Generalized Gross spaces exist in ZFC}

\bigskip \bigskip

Here we consider a natural generalization of the Gross property and will
obtain results in ZFC. 

Let $F$ be a field of arbitrary cardinality and $\l E, \Phi \r $ a
symmetric bilinear space over $F$, and let $\lambda $ be an infinite cardinal. We say that $\l E, \Phi \r $ is a
$\lambda $-{\it Gross space} if it has the following property:
\bigskip

$(\ast \ast \ast )$  for all subspaces $U\subseteq E$ of dimension
$\geq \lambda $: dim$U^\perp <$ dim$E$

\bigskip
So a Gross space is a $\o $-Gross space. By Lemma 1 in $\S 3$, $|F|^\la
$ is an upper bound for the dimension of a $\la $-Gross space over
$F$.

By results from [B/G] and [B/Sp] the existence of Gross spaces of
dimension $\aleph _1$ over
countable or finite
fields is independent of ZFC. We will prove here that in ZFC an
$\o _1$-Gross space of dimension $\aleph _1$ can be constructed over any
countably infinite
field. 

We need the following fancy Lemma:
\bigskip

{\bf Lemma 1.}  {\it Let $F$ be an infinite field. There exists a
sequence $\langle a_n :n<\o \rangle $ of elements of $F$ such that,
whenever $\langle k_{ij}: 1\leq i,j \leq m\rangle $ is a finite sequence of
pairwise distinct integers, then

$$\left |\matrix{ a_{k_{11}} & \ldots & a_{k_{1m}} \cr
                  \vdots     &        & \vdots     \cr
                  a_{k_{m1}} & \ldots & a_{k_{mm}} \cr }\right | \ne
0$$}

\bigskip

{\bf Proof.}  By a fusion argument. Start with a sequence $\langle
a^1_n :n<\o \rangle $ of pairwise distinct scalars in $F\setminus \{
0\} $. Then clearly every subsequence satisfies the conclusion of the
Lemma for $m=1$. 

Now assume $\langle a^m_n:n<\o \rangle $ has been constructed such
that the conclusion of the Lemma holds for any $1\leq m'\leq m$. Let $n_i=i$
for any $i\leq (m+1)^2-1$. Assume $n_0, \ldots ,n_l$ have been chosen.

Let $\langle i_0,j_0\rangle \in \{ 1,\ldots ,m+1\} ^2 $ and let
$\langle k_{ij}: \langle i,j\rangle \in \{ 1,\dots ,m+1\} ^2 \setminus
\langle i_0 ,j_0 \rangle \rangle $ be a sequence of pairwise distinct
integers in $\{ n_0,\ldots ,n_l\} $. By induction hypothesis, the
equation

$$\vcenter{\offinterlineskip
\def\strt{\vrule width 0pt height 10pt depth 6pt}
\let\\\displaystyle
\ialign{\hfil\ $\\#$\hfil\ &
	\vrule#&\hfil\ \ $\\#$\hfil\ &
			\hfil\ $\\#$\hfil\ &
				\hfil\ $\\#$\hfil\ \  &
			\vrule#&\qquad\strt$\\#$\cr
         &\omit&        &j_0\cr
         &\omit&        &\downarrow\strut\cr
         && a^m_{k_{11}}      &\cdots& a^m_{k_{1m+1}}    && \cr
i_0 \ \to\ && \vdots & X     & \vdots && = \ 0 \cr
         && a^m_{k_{m+11}}    &\cdots& a^m_{k_{m+1m+1}}    && \cr}}
 $$

\noindent in one variable $X$ which has coordinates $\langle i_0
j_0\rangle $ has exactly one solution. (Expand the determinant by the
$i_0$th row. Then the cofactors are all nonzero.)

Now choose $n_{l+1}>n_l$ minimal such that for any $n\geq n_{l+1}$, $a^m_n$ is distinct from
all the finitely many solutions of the equation above obtained by
running through all possible $\langle i_0,j_0\rangle $'s and $\langle
k_{ij}: \langle i,j \rangle \in \{ 1,\ldots m+1\} ^2\setminus \langle
i_0 j_0 \rangle \rangle $'s.

Now define $a^{m+1}_l=a^m_{n_l}$, $l<\o $.

Finally define

$$a_n =a^m_n \hbox{ iff }n\in \{ m^2-1,\ldots ,(m+1)^2-2\} $$

It is easily checked that $\langle a_n :n<\o \rangle $ is as
desired.\hfill $\square $
\bigskip

For uncountable fields an analog of Lemma 1 is true which is more
obvious. This is essentially [B/Sp, Lemma 2, $\S 4$].

\bigskip

{\bf Lemma 2.}  {\it Let $F$ be a field of uncountable cardinality
$\la $. There exists a sequence $\l a_\nu :\nu <\la \r $ of elements
of $F$ such that such that,
whenever $\langle \nu_{ij}: 1\leq i,j \leq m\rangle $ is a finite sequence of
pairwise distinct ordinals $<\la $, then

$$\left |\matrix{ a_{\nu _{11}} & \ldots & a_{\nu _{1m}} \cr
                  \vdots     &        & \vdots     \cr
                  a_{\nu _{m1}} & \ldots & a_{\nu _{mm}} \cr }\right | \ne
0$$}

\bigskip

{\bf Proof.}  Let $\l a_\nu :\nu <\la \r $ be a one-to-one enumeration
of a transcendence base of $F$ over its prime field. Prove by
induction on the matrix size $m$ that the conclusion is satisfied.
Expand the determinant by, say, the first row. Then by induction
hypothesis the cofactors are all nonzero, and the transcendentals of
the first row do not occur in them.

\hfill $\square $

\bigskip

{\bf Theorem 1.}  {\it A $\o _1$-Gross space of dimension $\aleph _1$ exists
over every countably infinite field.}

\bigskip

{\bf Proof.}  Let $F$ be a countably infinite field. Choose scalars $\langle a_n:n<\o \rangle $ in $F$ as in
Lemma 1. Let $\langle A_\alpha :\alpha <\omega _1 \rangle $ be an almost
disjoint family of subsets of ran$\langle a_n :n<\omega \rangle $.
For each $\alpha <\o _1$ choose a one-to-one function $f_\alpha :\al 
\rightarrow A_\alpha $. Let $E$ be a vector space over $F$ of
dimension $\aleph _1$ and let $\langle e_\alpha :\alpha <\o _1\rangle $
be a basis of $E$. 

On $E$, define a symmetric bilinear form as follows:
For $\alpha <\beta <\o _1$ define

$$\Phi (e_\alpha ,e_\beta )=f_\beta (\alpha )$$

$ \Phi (e_\alpha 
,e_\alpha )$ may be defined arbitrarily.

We claim that $\l E,\Phi \r $ is as desired. Using a $\Delta
$-system argument and the bilinearity of $\Phi $, it is not difficult
to see that it is enough to prove the following:

\bigskip

($\star $) Assume that $U$ is a subspace spanned by a basis $\langle
y_k :k<\omega \rangle $ such that the sets in $\langle$supp$(y_k):k<\o
\rangle $ are pairwise disjoint and of the same cardinality, then
dim$U^\perp <\omega _1$. 

\bigskip

By way of contradiction assume that $U$ and $\langle y_k:k<\omega
\rangle $ are as in $(\star )$ but $U^\perp $ contains $\o _1$ many
linearly independent vectors $\langle z_\iota :\iota <\o _1\rangle $.

Let 

$$y_k=\sum_{l=1}^m b_{kl}e_{\alpha (k,l)}$$

We may assume that each $z_\iota $ has the same nonzero coefficients
in its representation, say 

$$z_\iota =\sum_{l=1}^n c_l e_{\beta (\iota, l)}$$

\noindent and for any $\iota <\nu <\o _1$ the sets supp$(z_\iota )$
and supp$(z_\nu )$ are disjoint and supp$(z_\iota )$ is disjoint from
$\sup\bigcup_{k<\o
}$supp$(y_k)$.

We claim that not even the first $m\cdot n$ many $z_\iota $'s are in
$U^\perp $. In order to see this, let us compute

$$\vbox{\halign{# & # & # & \hfill # \hfill & # & \hfill # \hfill & #
\cr 

$\Phi (z_\iota ,y_k)$&$=$&&$c_1b_{k1}f_{\beta (\iota ,1)}(\alpha
(k,1))$& 

$+\ldots +$&$c_1b_{km}f_{\beta (\iota ,1)}(\alpha (k,m))$&\cr

&&&$\vdots $&&$\vdots $&\cr

&&$+$&$c_nb_{k1}f_{\beta (\iota ,n)}(\alpha (k,1))$&$+\ldots +$&$c_nb_{km}f_{\beta (\iota ,n)}(\alpha (k,m))$&\cr}}$$

Since the sets $A_\alpha $ are almost disjoint, the sets
supp$(y_k)= \{ \alpha (k,1),\ldots ,\alpha (k,m)\} $ are pairwise
disjoint, the functions $f_\alpha $ are one-to-one and the scalars
$\langle a_n:n<\o \rangle $ are chosen satisfying the conclusion of
Lemma 1, there is $k<\o $ such that for $\iota <m\cdot n$ the sets 

$$\{ f_{\beta (\iota ,i)}(\alpha (k,j)): 1\leq i\leq n, \, 1\leq j
\leq m\} $$

\noindent are pairwise disjoint, contain $m\cdot n$ elements and

$$\left\vert\matrix{ f_{\beta(0,1)}(\alpha (k,1)) & \ldots &f_{\beta
(0,n)}(\alpha (k,m)) \cr
\vdots & &\vdots \cr
f_{\beta (mn-1,1)}(\alpha (k,1)) & \ldots & f_{\beta (mn-1,n)}(\alpha
(k,m)) \cr}
\right\vert \ne 0 \eqno(1)$$

But then $y_k$ is not orthogonal on every $z_\iota $, $\iota <m\cdot
n$, since otherwise the vector

$$\langle c_1b_{k1}, \ldots , c_nb_{km}\rangle \in F^{mn}$$

\noindent would be a nontrivial solution of the system of homogeneous
equations $Ax=0$ where $A$ is the $mn \times mn-$matrix in (1).\hfill
$\square $ 

\bigskip

The upper bound for the dimension of a $\o _1$-Gross space
over a countable field is $\o ^{\o _1}=2^{\o _1}$, by Lemma 1, $\S 3$.
Hence, the largest dimension in which such a space conceivably can be
constructed in ZFC is $\o _2$. Theorem 6 below will show that in fact
such a construction is possible.

For the following question we have no answer:

\bigskip

{\bf Question.} Does there exist a $\o _1$-Gross space of uncountable
dimension over any finite field?

\bigskip

Next we will show that this question has a positive answer for regular
cardinals larger than $\o _1$. For this, recent work of the first author on
colouring pairs of ordinals will be applied. Let $\la ,\mu ,\ka ,
\theta $ be cardinals such that $\la $ is infinite and $\la \geq \mu
\geq \ka +\theta $. Let Pr$_0(\la ,\mu ,\ka ,\theta )$ be the
following statement:

\bigskip
Pr$_0(\la ,\mu ,\ka ,\theta )$:  There exists a function $c:[\la
]^2\rightarrow \ka $ such that, if $\xi <\theta $ and for $\nu <\mu $,
$\l \al _{\nu,\zeta }:\zeta <\xi\r $ is a strictly increasing sequence
of ordinals $<\la $, the $\al _{\nu ,\zeta }$ distinct, and $h: \xi \times \xi\rightarrow  \ka $, then there are $\nu _1<\nu _2<\mu $ such that
for all $\zeta _1,\zeta _2 <\xi $ we have $c\{ \al _{\nu _1, \zeta
_1}, \al _{\nu _2, \zeta _2}\}=h\l \zeta _1, \zeta _2\r $. 

\bigskip 

The following theorem, which is due to the first author, states that
Pr$_0(\la ^+, \la ^+,\la ^+, \o )$ holds for all regular uncountable cardinals.
For regular $\la >\o _1 $ this is proved in [Sh4, Corollary (a), 
Section 4, p.100]; see Definition 1, p.95. For $\la =\o _1$ this was 
proved in [Sh7, Theorem 1.1].

\bigskip

{\bf Theorem 2.}  Pr$_0(\la ^+, \la ^+,\la ^+, \o )$ {\it and
hence} Pr$_0(\la ^+, \la ^+, 2, \o )$ {\it holds for every regular
uncountable cardinal $\la $.}

\bigskip
For more on Theorem 2 (for $\la >\o _1$) see also [Sh5]. 
It appears (when we specify $\theta =
\o $) in [Sh5, III, 4.8(1), p.177]. Moreover, 4.8(2) gives similar
results for inaccessibles with stationary subsets not reflecting in 
inaccessibles. 

\bigskip

{\bf Theorem 3.}  {\it Let $\la $ be a regular uncountable cardinal. A
$\la
^+$-Gross space of dimension $\la ^+$ exists over any field.}

\bigskip

{\bf Proof.}  Let $F$ be an arbitrary field of size at most $\la $
(for larger fields see Theorem 4 below) and
$E$ a vector space
over $F$ of dimension $\la ^+$. Choose a basis $\l e_\al :\al <\la
^+\r $. We define a symmetric bilinear form on $E$ using the colouring
given by Pr$_0(\la ^+,\la ^+,2,\o )$. For $\{ \al ,\be \} \in [\la
^+]^2$ we set

$$\Phi (e_\al ,e_\be )=c\{ \al ,\be \} $$

The angles $\Phi (e_\al ,e_\al )$ may be defined arbitrarily. 

Assume that $\l E,\Phi \r $ is not $\la ^+$-Gross. So we may find
families $\l y_\nu :\nu <\la ^+\r $ and $\l z_\nu :\nu <\la ^+\r $ of
linearly independent vectors such that $\Phi (y_\nu ,z_{\nu '})=0$
always. Without loss of generality we may assume that $\l $dom$(y_\nu
): \nu <\la ^+\r $ and $\l $dom$(z_\nu ):\nu <\la ^+\r $ are families
of pairwise disjoint sets, the sets in the first one all of size $n$
and those in the latter all of size $m$, and if $\nu <\nu '$ and $\al
\in $ dom$(y_\nu )$, $\be \in $ dom$(z_\nu )$, $\gamma \in $
dom$(y_{\nu '})$ then $\al <\be <\gamma $. Furthermore we may assume
that each $y_\nu $ and each $z_\nu $ has the same nonzero coefficients
in its representation, say

$$y_\nu =\sum_{l<n} b_{l}e_{\beta (\nu ,l)}$$ 

$$z_\nu =\sum_{l<m} c_l e_{\gamma (\nu, l)}$$

\noindent We assume that $\beta (\nu ,l)$ and $\gamma (\nu ,l)$
increase with $l$. Let $\l \al _{\nu ,i}:i<n+m\r $ be the
increasing enumeration of dom$(y_\nu )\cup $ dom$(z_\nu )$. 

It is easy to define $h:(n+m) \times (n+m)\rightarrow 2$ such that

$$\sum_{l< n}\sum_{l'< m}b_l c_{l'}h\l l, n+l'\r \ne
0$$

But now by Pr$_0(\la ^+,\la ^+,2,\o )$ we may find $\nu _1<\nu _2<\la
^+$ such that for all $i,j<n+m$ we have $c\{ \al _{\nu _1,i}, \al
_{\nu _2,j}\} =h\l i,j\r $. We conclude $\Phi (y_{\nu _1},z_{\nu
_2})\ne 0$, a contradiction.\hfill $\square $

\bigskip

For uncountable fields of cardinality equal to the dimension of the space a
much stronger version of Theorems 1 and 3 was proved
in [G/O]:

\bigskip

{\bf Theorem 4 [G/O, Thm. 1, $\S 1$].}  {\it Let $F$ be a field of
uncountable cardinality $\la $. A strong Gross space of dimension $\la
$ over $F$ exists.}

\bigskip

Here a strong Gross space is a space $\l E,\Phi \r $ satisfying the
following property:

\bigskip 

$(\ast \ast )$  for all subspaces $U\subseteq E$ of infinite
dimension: dim$U^\perp \leq \aleph _0$

\bigskip

For more on strong Gross spaces see [Sp4]. 

Constructing a Gross
space gets the more difficult, the larger its dimension is compared with
the size of its base field. For $\lambda =\o $, the existence of a
$\lambda 
$-Gross space of dimension $\la ^+$ over some field of size $\lambda $
is independent of ZFC. The next theorem shows that if $\la $ is
uncountable, then such a space can be constructed in ZFC. If
cf$(\lambda )  >\o $, a similar construction as in
Theorem 1 is possible. If cf$(\lambda )=\o $, the almost disjoint sets
used for the construction must be chosen more carefully. For this we
will need the following
result from the first author's work on cardinal arithmetic. It is obtained
from [Sh6, 6.2] with the ideal of bounded subsets of $\kappa $ and the
$f_\alpha $'s coming from [Sh5, II, 1.5, p.50].

\def \la {\lambda } \def \i {\iota } \def \ka {\kappa }
\bigskip

{\bf Lemma 3 [Sh6, 6.2].}  {\it Assume that $\la $ is a singular
cardinal and} cf$(\la ) =\kappa $. {\it There exists an increasing sequence
$\l \la _\i : \i <\ka \r $ of regular cardinals $\la _\i >\ka $ with
limit $\la $ and a family $\l f_\al :\al <\la ^+\r $ of functions in $
\prod_{\i <\ka }\la _\i $ such that for all
$\al < \beta < \la ^+$:

(1) $f_\al <f_\beta $ modulo bounded subsets of
$\ka $, and

(2) if $\l u_\zeta :\zeta <\la ^+\r $ is a sequence of pairwise
disjoint nonempty subsets of $\la ^+ $, each $u _\zeta $ of size
$<\ka $, and $\mu <\la $ is a cardinal, then there exists $B\subseteq \la
^+$ of size $\mu $ and $\i ^* <\ka $ such that for all $\zeta , \, \xi
\in B$ and $\i <\ka $:

\itemitem{(i)} if $\zeta <\xi $, then $\sup u_\zeta <\min u_\xi $, and

\itemitem{(ii)} if $\al ,\, \beta \in \bigcup_{\zeta \in
B}u_\zeta $, $\al <\beta $ and $\i ^* \leq \i <\ka $, then $f_\al (\i )<f_\beta (\i )$.}

\def \be {\beta }
\bigskip

{\bf Theorem 5.}  {\it Let $\lambda $ be an uncountable cardinal and
$F$ a field of size $\lambda $. Then a $\lambda $-Gross space over $F$ of
dimension $\lambda ^+$ exists.}
\def \z {\zeta } \def \la {\lambda } \def \i {\iota } 
\bigskip

{\bf Proof.} Let $F$ be a field of cardinality $\la $. Let $\l
a_\nu :\nu <\lambda \r $ be a transcendence base of $F$ over its
prime field $F_0$. Hence, $F$ is an algebraic extension of $F_0(\l
a_\nu :\nu <\la \r ).$ Let $E$ be a $F$-vector space of dimension $\la ^+$, spanned by a basis
$\l e_\al :\al <\la ^+\r $.

\bigskip

{\it Case 1:}  Assume cf$(\la )>\o $. Choose $\l A_\al : \al <\la ^+\r
$, a family of almost disjoint subsets of $\la $. I.e. each $A_\al $
has size $\la $ and if $\al ,\, \beta $ are distinct then $|A_\al \cap
A_\beta |<\la $ (see [J, p.252]). For each $\al <\la ^+ $ choose a
one-to-one function $g_\al :\al \rightarrow A_\al $. 

On $E$, define a symmetric bilinear form as follows: For $\al <\beta
<\la ^+$ let 

$$\Phi (e_\al ,e_\beta )= a_\nu \hbox{ if and only if }g_\beta (\al
)=\nu $$ 

The angles $\Phi (e_\al ,e_\al )$ may be defined arbitrarily. 

We claim that $\l E,\Phi \r $ is a $\la $-Gross space. Since cf$(\la
)>\o $, in every subspace $U\subseteq E$ of dimension $\la $
we may find $\la $ many
vectors such that the supports of each of them has the same size. We
also may assume that these supports are pairwise disjoint (see Case 2).
Hence it is enough to prove the following:

\bigskip

$(\star \star )$ Assume that $U$ is a subspace spanned by a basis $\langle
y_\nu :\nu <\la \rangle $ such that the sets in $\langle$supp$(y_\nu
):\nu <\la 
\rangle $ are pairwise disjoint and of the same cardinality. Then
dim$U^\perp <$ dim$E$. 

\bigskip

The proof of this is completely
analogous to that of $(\star )$ in the proof of Theorem 1, if we
use Lemma 2 instead of Lemma 1.

\bigskip

{\it Case 2:}  Assume cf$(\la )=\o $. In this case, we have to choose
the almost disjoint sets $\l A_\al :\al <\la ^+\r $ more carefully.
Let $\l \la _n :n<\o \r $ and $\l f_\al :\al <\la ^+\r $ be as in
Lemma 3 where $\kappa =\o $.

Let $\l \mu _n :n <\o \r $ be increasing,
continuous and with limit
$\la $, such that $\mu _n <\la _n $ for every $n $. For each $\al
<\la ^+$ define:

$$A_\al =\bigcup_{n <\o }[\mu _n +\mu _n \cdot f_\al (n ), \,
\mu _n +\mu _n \cdot (f_\al (n )+1) ) $$

\noindent and let $g_\al :\al \rightarrow A_\al $ be one-to-one and onto.

Now let $E$ be a vector space over $F$ of dimension $\la ^+$, spanned
by a basis $\l e_\al :\al <\la ^+\r $. On $E$, define a symmetric
bilinear form as follows: If $\al <\beta <\la ^+$, let:

$$\Phi (e_\al ,e_\beta )=a_{g_\beta (\al )}$$

\noindent The angles $\Phi (e_\alpha ,e_\alpha )$ may be defined
arbitrarily. 

We claim that $\l E, \Phi \r $ is a $\la $-Gross space. Assume that
this is not true. So there is a subspace $U$, spanned by a basis $\l
y_\nu : \nu <\la \r $, such that $U^\perp $ contains a family $\l z_\z
: \z <\la ^+\r $ of linearly independent vectors. Without loss of
generality we may and do assume that $\l $supp$(y_\nu ):\nu <\la \r $ is
a family of pairwise disjoint sets. This may be seen as follows: For
each $n <\o $ there exists $Y_n \subseteq [\la _n , \la _n ^+
) $ of size $\la _n ^+$ such that $\l $supp$( y_\nu) :\nu \in Y_n \r $
is a $\Delta $-system. By forming linear combinations, we find $X_n
\subseteq $ span$\l y_\nu :\nu \in Y_n \r $ such that $\l $supp$(x):
x\in X_n \r $ is a pairwise disjoint
family. By cutting off from each $X_n $ at most $\sum_{n '<n }\la
_{n '}$ many elements, we obtain a family as desired.

Let

$$y_\nu =\sum_{l=1}^{m_\nu }b_{\nu l}e_{\alpha (\nu ,l)}$$

\noindent where each $b_{\nu l}$ is nonzero. Choose $\al ^*<\la ^+$
such that $\forall \nu <\la \forall 1\leq l\leq m_\nu $ $\al (\nu ,l)<\al ^*$.

Without loss of generality we may assume that there are nonzero $c_1,
\ldots , c_n$ such that each $z_\z $ has a representation 

$$z_\z =\sum_{l=1}^nc_l e_{\beta (\z ,l)}$$

\noindent and for all $\z <\eta $ we have $\al ^*<\beta (\z ,
1)<\ldots <\beta (\z, n)<\beta (\eta ,0)<\ldots <\beta (\eta , n)$.

By applying Lemma 3 to the family $\l $supp$(z_\z ): \z < \la ^+ \r $
and $\mu =\o $, we find $B\in [\la ^+]^\o $ and $n^*<\o $ such that
(ii) from the conclusion of Lemma 3 holds.

Note that now $\l A_\beta \setminus \la _{n^*}: \beta \in
\bigcup_{\zeta \in B}$dom$(z_\zeta )\r $ is a family of disjoint sets.
Since also $\l $dom$(y_\nu ):\nu <\la \r $ is a disjoint family and
the $g_\al $'s are one-to-one we may find $y_\nu $ such that for all
$\beta \in \bigcup_{\zeta \in B}$dom$(z_\zeta )$ we have $g_\beta
($dom$(y_\nu ))\subseteq A_\beta \setminus \la _{n^*}$. 

But now not even $m_\nu \cdot n$ many vectors $z_{\z _1},\ldots ,
z_{\z _{m_\nu n}}$, where each $\z _i\in B$, are orthogonal on $y_\nu $.
For if this would be true, then $\l
b_1c_1,\ldots ,b_{m_\nu }c_n\r \in F^{m_\nu n}$ would be a nontrivial
solution of
the equation $Ax=0$ where $A$ is the $(m_\nu n\times m_\nu n)$-matrix

$$\left(\matrix{ \Phi (e_{\al (\nu ,1)}, e_{\be (\z _1,1)}) & \ldots &
\Phi (e_{\al (\nu ,m_\nu )},e_{\be (\z _1,n)}) \cr 
\vdots & & \vdots \cr
\Phi (e_{\al (\nu ,1)} ,e_{\be (\z _{m_\nu n},1)}) & \ldots & \Phi (e_{\al
(\nu ,m_\nu )}, e_{\be (\z _{m_\nu n},n)})\cr }\right) $$

But by construction and Lemma 2, this matrix has nonzero determinant,
a contradiction.\hfill $\square $

\bigskip

By Lemma 1 in $\S 3$, the upper bound for a $\la $-Gross space over
$F$ is $|F|^\la $. Hence in ZFC, the dimension of a space as in
Theorem 5 cannot be enlarged. But the following theorem shows that it
can be enlarged at the loss that we only obtain a $\la ^+$-Gross
space. An upper bound for a $\la ^+$-Gross
space over a field of size $\la $ is $\la ^{\la ^+}=2^{\la ^+}$. So
again, in ZFC we cannot expect a better result. Theorem 6 is true also
for countable fields, so we get a complement to Theorem 1.

The construction will use filtrations given by the following Lemma:
\def \g {\gamma }
\bigskip

{\bf Lemma 4.}  {\it Let $\la $ be a cardinal.
\bigskip

(A) There exists a family $\l A^\al _\g :\g <\la ^+, \al <\la ^{++}\r
$ of sets of size $\leq \la $ such that for any $\al , \beta
<\la ^{++}$ and $\gamma ,\delta < \la ^{+}$ the following requirements are
satisfied: 

(1) $\l A^\al _\gamma :\gamma <\la ^+ \r $ is increasing,
continuous such that $\al =\bigcup_{\g <\la ^+ }A^\al _\g $ and
$A^\al _0=\emptyset $.

(2) if $\g <\delta $ then $A^{\al } _{\gamma } \subseteq A^{\al }_{\delta }$

(3) if $\al <\beta $ and $\al \in A^\be _\g $ then $A^{\al }_{\g } =
A^{\be }_{\g } \cap \al $.

\bigskip

(B) Let $\l C_\al :\al <\la ^{++}\r $ be a family of clubs of $\la
^{+}$. There exists a family $\l C_\al ':\al <\la ^{++}\r $ of clubs
and a sequence $\l A^\al _\g :\g \in C_\al ' , \al <\la ^{++}\r
$ of sets of size $\leq \la $ such that for any $\al , \beta
<\la ^{++}$ and $\gamma ,\delta < \la ^{+}$ we have $C_\al '\subseteq C_\al
$, and if $\al <\be $ the set $C_\be '\setminus C_\al '$ is bounded and the
following requirements are
satisfied: 

(1) $\l A^\al _\gamma :\gamma \in C_\al '\r $ is increasing,
continuous such that $\al =\bigcup_{\g \in C_\al '}A^\al _\g $ and
$A^\al _0=\emptyset $.

(2) if $\gamma ,\delta \in C_\al ' $ and $\gamma <\delta $ then $A^{\al }_{\gamma }
\subseteq A^{\al }_{\delta } $

(3) if $\al <\beta $ and $\al \in A^{\beta }_{\gamma }$ then $C_\be '\setminus
C_\al ' \subseteq \g $ and $A^{\al }_{\gamma } =A^{\beta }_{\gamma }
\cap  \al $.} 

\bigskip

{\bf Proof.}  We prove only (A). The proof of (B) is similar.

Assume that $\l A^\nu _\g : \g <\la ^+\r $ for $\nu <\al $ have been constructed
satisfying the requirements. Let $\l \al _\zeta :\zeta <\la ^+\r $ be
an enumeration of $\al $. Then the set 

$$\{ \g <\la ^+: (\forall \zeta _1, \zeta _2 <\la ^+)\hbox{ if }
\zeta _1<\zeta _2 \hbox{ and }\al _{\zeta _1}<\al _{\zeta _2} \hbox{
then } \al _{\zeta _1}\in A_\g ^{\al _{\zeta _2}}\} $$

\noindent
clearly is a club of $\la ^+$. Let $\l \g _\varepsilon
:\varepsilon <\la ^+\r $ be its increasing enumeration and define
$A_0^\al =\emptyset $, and 

$$A^\al _\g =\bigcup_{\zeta <\g _\varepsilon }A_\g ^{\al _\zeta }$$

\noindent if and only if $\g _\varepsilon <\g \leq \g
_{\varepsilon +1}$.   

It is not difficult to see that this works.\hfill $\square $

\bigskip

{\bf Theorem 6.}  {\it Let $\la $ be an infinite cardinal. A $\la ^+
$-Gross space of dimension $\la ^{++}$ exists over every field of
size $\la $.}

\def \be {\beta } \def \i {\iota } \def \g {\gamma } \def \d {\delta }
\bigskip

{\bf Proof.}  The proof is divided into two cases:

\bigskip

{\it Case 1:}  There exists a family $\l C_\al :\al <\la ^{++}\r $ of
clubs of $\la ^+$ such that for every club $C$ of $\la ^+$, for some
$\al $ and $\la ^+$ many $\g \in C$ we have $\min (C\setminus (\g
+1))\leq \min (C_\al \setminus (\g +1))$.  Notice that this is
equivalent to the assumption that there exist $\la ^{++}$ clubs of
$\la ^+$ such that no club gets eventually inside all of them.

By Lemma 4(B), without loss of generality we may assume that for $\al <\be $ the set $C_\be
\setminus C_\al $ is bounded in $\la ^+$ and we have families
$\l A^\al _\g :\g \in C_\al \r $ satisfying (1), (2) and (3). 

Now let $F$ be a field of size $\la $ and $E$ a vector space over $F$,
spanned by a basis $\l e_\al :\al <\la ^{++}\r $. 

On $E$, we will define a symmetric bilinear form $\Phi $ by induction
on $\al <\la ^{++}$, using the filtrations above. Assume that $\Phi (e_\beta ,e_{\beta '})$ has
been defined for $\beta \leq \beta '<\al $ and $\Phi (e_\beta ,e_\al
)$ has been defined for $\beta \in A^\al _\g $.

Let $\d $ be the successor of $\g $ in $C_\al $. If $A^\al _\d
\setminus A^\al _\g $ has size less than $\la $, then $\Phi (e_\beta
,e_\alpha )$, for $\beta \in  A^\al _\d
\setminus A^\al _\g $, and $\Phi (e_\al ,e_\al )$ may be defined
arbitrarily. Otherwise, choose one-to-one enumerations such that
$F=\hbox{ ran}\l a_\nu :\nu <\la \r $, $A^\al _\g =\hbox{ ran}\l \zeta _\nu
:\nu <\la \r $ and $A^\al _\d \setminus A^\al _\g =$ ran$\l \xi
_\nu :\nu <\la \r $.

Assume that $\Phi (e_{\xi _{\nu '}}, e_\al )$ has been defined for $\nu
'<\nu $. 

Define $\Phi (e_{\xi _\nu },e_\al )$ such that $\Phi (z,y)\ne 0$
whenever $z$ and $y$ are vectors such that $\al \in $ dom$(z)\subseteq \{ \al \} \cup \{ \zeta _{\nu '} :\nu ' <\nu \} $ and
$\xi _\nu \in $ dom$(y)\subseteq \{ \xi _{\nu '}:\nu '\leq \nu \}
$ and all coefficients of $z$ and $y$ are contained in $\{ a_{\nu '}:
\nu '<\nu \} $. 
There are less than $\la $ many such $z$ and $y$. Hence, this
definition can be fulfilled.

We will show that $\l E, \Phi \r $ is a $\la ^+$-Gross space. Assume
that this is not true. So there is a subspace $U\subseteq E$ of
dimension $\la ^+$, say $U=\hbox{ span}\l y_\nu :\nu <\la ^+\r $, such
that in $U^\perp $, there is a family of $\la ^{++}$ linearly
independent vectors, say $\l z_\i :\i <\la ^{++}\r $. 

Without loss of generality, we may assume that $\l \hbox{supp}(y_\nu
):\nu <\la ^+\r $ is a family of pairwise disjoint sets of the same
cardinality and each $y_\nu $ has the same nonzero coefficients in its
representation. Mutatis mutandis, we may assume the same for the $z_\i
$'s. 

Choose $\al ^*<\la ^{++}$ such that $U\subseteq \hbox{ span}\l e_\al : \al <\al
^*\r $.

We may certainly find a club $C\subseteq C_{\al ^*}$ such that, if
$\gamma <\d $ are successive members of $C$, then 

$$|\{ \nu <\la ^+ :
\hbox{ dom}(y_\nu ) \subseteq A^{\al ^*}_\d \setminus A^{\al ^*}_\g \}
| = \la \eqno(1)$$

Since we are in Case 1 there exists $\al _C
\in (\al ^*, \la ^{++})$ such that for $\la ^+$ many $\g \in C$ we have
$ \min (C\setminus (\g +1))\leq \min (C_{\al _C } \setminus (\g +1)) $.
       
Now let $\i <\la ^{++}$ such that dom$(z_\i ) \not\subseteq \al _C$
and let $\al $ the largest member of dom$(z_\i )$. Hence $\al \geq \al
_C$, $C_\al \setminus C_{\al _C}$ is bounded, and hence 
$S=\{ \g \in  C :  \min
(C\setminus (\g +1))\leq \min (C_\al \setminus (\g +1))\} $ has size
$\la ^+$. 

Now choose $\g '\in S $ so large that, if $\g $ is the maximal element
of $C_\al $ such that $\g \leq \g '$, then   

$$\{ \al ^*\}\cup \hbox{ dom}(z_\i )\setminus \{ \al \} \subseteq
A^\al _{\g } \eqno(2)$$

\noindent Let $\d $ be the successor of $\g $ in $C_\al $. So
clearly $\min (C\setminus (\g '+1))\leq \d $. 

By construction of the filtrations $\l A^\al _\nu :\nu \in C_\al \r $
and by (1) and (2) we know that

$$|\{ \nu <\la ^+: \hbox{ dom}(y_\nu ) \subseteq A^\al _\d \setminus
A^\al _\g \} |=\la \eqno(3)$$

Now look at the definition of $\Phi (.,e_\al )$ on $A^\al _\d \setminus
A^\al _\g $. Choose $\nu <\la $ and $\mu <\la ^+$ such that 

(i) the coefficients of $z_\i $ and those of $y_\mu $ (and hence of
any $y_{\mu '}$) are contained in
$\{ a_{\nu '}:\nu '<\nu \} $,

(ii) dom$(z_\i )\setminus \{ \al \} \subseteq \{ \zeta _{\nu '}:\nu
'<\nu \} $, and

(iii) $\xi _\nu \in \hbox{ dom}(y_\mu )$ and$\hbox{ dom}(y_\mu )\setminus \{ \xi _\nu \}
\subseteq \{ \xi _{\nu '}:\nu '<\nu \} $.

By (3), such a choice is possible. But now $\Phi (e_{\xi _\nu } ,e_\al
) $ was defined such
that $\Phi (z_\i ,y_\mu )\ne 0$, a contradiction.    

\bigskip

{\it Case 2:}  There exists no family of clubs as in Case 1.

Let $\l A^\al _\g :\g <\la ^+ , \al <\la ^{++}\r $ be as in Lemma 4(A). 

Let $f_\al :\la ^+ \rightarrow \la ^+ $ be defined by $f_\al (\g )=$
o.t.$(A^\al _{\g +1} )$. Note that if $\al <\be $, then $f_\al (\g
)<f_\be (\g )$ for every $\g $ so large that $\al \in A^\be _\g $. 

Let $C_\al =\{ \g <\la ^+: \g $ is a limit ordinal and $(\forall \d
<\g ) f_\al (\d )<\g \} $. Clearly $C_\al $ is a club.

Since $\l C_\al :\al <\la ^{++}\r $ cannot serve for Case 1, there
exists a club $C$ of $\la ^+$ such that 

$$\forall \al <\la ^{++}\exists \g <\la ^+\forall \d \in (\g ,\la ^+)\cap C
\, \, [\min (C\setminus (\d +1))>\min (C_\al \setminus (\d +1))]
\eqno(4)$$

\def \e {\varepsilon }
Let $\l \g _\e :\e <\la ^+\r $ be the increasing enumeration of $C\cup
\{ 0\} $. Now define $f^*:\la ^+\rightarrow \la ^+$ by:

$$f^*(\g )=\g _{\e +2} \hbox{ if and only if }\g _\e \leq \g <\g _{\e
+1}$$

Using the definition of $C_\al $ and (4), it is easy to see that:

$$\forall \al <\la ^{++}\exists \g <\la ^+ \forall \d \in (\g , \la
^+) \, \, [f_\al (\g )<f^*(\g )]\eqno(5)$$

Now let $F$ be a field of size $\la $ and $E$ a vector space over $F$
of dimension $\la ^{++}$, spanned by a basis $\l e_\al :\al <\la
^{++}\r $. 

Let $\l F_\g :\g <\la ^+\r $ be an almost disjoint family of subsets
of $F$ (so each $F_\g $ has size $\la $ and $F_\g \cap F_\d $ has size
$<\la $ for any $\g \ne \d $) such that:

{\it in case} $\la $ is uncountable, $\bigcup F_\g $ is a set of
elements which are algebraically independent over the prime field of
$F$

{\it in case} $\la =\o $, $\bigcup F_\g $ is a subset of ran$\l a_n:
n<\o \r $ where $\l a_n:n<\o \r $ is a sequence as in Lemma 1.

For each $\g <\la ^+$, let $\l a^\g _{\e , \zeta }:\e <\zeta <f^*(\g
)\r $ be a one-to-one enumeration of $F_\g $.

In difference to the previous constructions in this paper, here the
angles $\Phi (e_\al ,e_\be )$ will be defined upwards, as follows:

For $\al <\la ^{++}$ and $\g <\la ^+$ define

$$W^\al _\g=\{ \be <\la ^{++}: \al \in A^\be _{\g +1}\setminus A^\be
_\g \} $$

\noindent So $\la ^{++}\setminus \al =\bigcup_{\g <\la ^+}W^\al _\g $, and this is
a disjoint union. 

Now let $\al <\be <\la ^{++} $. We define

$$\Phi (e_\al ,e_\be )=a^\g _{{o.t.}(A^\al _{\g +1}), {
o.t.}(A^\be _{\g +1})}$$

\noindent in case o.t.$(A^\be _{\g
+1}) <f^*(\g )$ where $\g $ is the uniquely determined ordinal such
that $\be \in W^\al _\g $. Note that o.t.$(A^\al _{\g +1})< $ o.t.$(A^\be
_{\g +1} )$ in this case. Otherwise, $\Phi (e_\al ,e_\be )$ may be defined
arbitrarily. Also $\Phi (e _\al ,e_\al )$ may be defined arbitrarily.

We claim that $\l E,\Phi \r $ is $\la ^+$-Gross. If this is not true,
then there is a subspace $U\subseteq E$ of
dimension $\la ^+$, say $U=$ span$\l y_\nu :\nu <\la ^+\r $, such
that in $U^\perp $, there is a family of $\la ^{++}$ linearly
independent vectors, say $\l z_\i :\i <\la ^{++}\r $. 

Again, we may assume that $\l $supp$(y_\nu
):\nu <\la ^+\r $ is a family of pairwise disjoint sets of the same
cardinality and each $y_\nu $ has the same nonzero coefficients in its
representation. Mutatis mutandis, we may assume the same for the $z_\i
$'s. We may also assume that each supp$(z_\iota )$ is disjoint from
$\al ^*$ where $\al ^*$ is chosen such $U\subseteq $ span$\l e_\al :\al <\al
^*\r $. Let   

$$y_\nu =\sum_{l=1}^m b_l e_{\al (\nu ,l)}$$

$$z_\i =\sum_{l=1}^n c_l e_{\be (\i ,l)}$$

By (5) and property A(3), Lemma 4, of the filtrations $\l A^\al _\g
:\g <\la ^+, \al <\la ^{++}\r $ we can choose $\g ^*<\la ^+$ such that

(j) for each $\be \in \bigcup_{\i
<\la }$dom$(z_\i )$ and $\g \in (\g ^*, \la ^+)$ $f_\be (\g )<f^*(\g
)$, and

(jj) if $\be _1, \, \be _2\in \bigcup_{\i
<\la }$dom$(z_\i ) $ and $\be _1<\be _2$, then $\be _1 \in A^{\be
_2}_{\g ^*} $ and $\al ^*\in A^{\be _1}_{\g ^*}$.

Now find $y_\nu $ such that dom$(y_\nu )\cap A^{\al ^*}_{\g ^*} =
\emptyset $. There are finitely many $\g _1,\ldots ,\g _p>\g ^*$ such
that 

$$\hbox{dom}(y_\nu )\subseteq \bigcup_{i=1}^p(A^{\al ^*}_{\g
_i+1}\setminus A^{\al ^*}_{\g _i})$$

Then for all $\be \in  \bigcup_{\i
<\la }$dom$(z_\i ) $ and $1\leq i\leq p$ we have 

$$\hbox{dom}(y_\nu )\cap (A^{\al ^*}_{\g
_i+1}\setminus A^{\al ^*}_{\g _i})=\hbox{dom}(y_\nu )\cap (A^{\be }_{\g
_i+1}\setminus A^{\be }_{\g _i})$$

\noindent Call this set $M_i$. So for $\al \in M_i $ and $\beta \in
 \bigcup_{\i
<\la }$dom$(z_\i ) $ we have $\be \in W^\al _{\g _i} $ and hence, by
(j) above, $\Phi
(e_\al ,e_\be )= a^{\g _i} _{{o.t.}(A^\al _{\g _i +1}), {
o.t.}(A^\be _{\g _i +1})} \in F_{\g _i}$. Notice that for distinct
$\al ,\al ' \in M_i $ we have o.t.$(A^\al _{\g _i+1})\ne $ o.t.$(A^{\al '}_{\g
_i+1})$, and for distinct $\be ,\be '\in  \bigcup_{\i
<\la }$dom$(z_\i ) $ we have o.t.$(A^\be _{\g _i+1})\ne $ o.t.$(A^{\be '}_{\g
_i+1})$. Hence, by almost-disjointness of the $F_{\g _i}$, we may find
$\i _1, \ldots ,\i _{m\cdot n}<\la $ such that for $1\leq j\leq
mn$ the sets

$$ \{  a^{\g _i} _{{o.t.}(A^\al _{\g _i +1}), {
o.t.}(A^{\be (\i _j,l)} _{\g _i +1})} : \al \in M_i, \, 1\leq i\leq p,
\, 1\leq l\leq n\} \eqno(6)$$

\noindent are pairwise disjoint and each of them contains $m\cdot n$
elements. If now 

$$\Phi (y_\nu ,z_{\i _j})= \sum_{1\leq l\leq m, 1\leq k\leq
n}b_lc_k\Phi (e_{\al (\nu ,l)}, e_{\be (\i _j,k)})=0$$

\noindent would be true for every $1\leq j\leq mn$, then $\l
b_1c_1,\ldots ,b_mc_n\r \in F^{mn}$ would be a nontrivial solution of
the equation $Ax=0$ where $A$ is the $(mn\times mn)$-matrix

$$\left(\matrix{ \Phi (e_{\al (\nu ,1)}, e_{\be (l_1,1)}) & \ldots &
\Phi (e_{\al (\nu ,m)},e_{\be (l_1,n)}) \cr 
\vdots & & \vdots \cr
\Phi (e_{\al (\nu ,1)} ,e_{\be (l_{mn},1)}) & \ldots & \Phi (e_{\al
(\nu ,m)}, e_{\be (l_{mn},n)})\cr }\right) $$ 

\noindent But the rows of this matrix are the sets in (6), hence its
determinant is nonzero, a contradiction.\hfill $\square $

\vfill\eject

\def\ltextindent#1{\hbox to \hangindent{#1\hss}\ignorespaces}
\def\litem{\par\noindent
               \hangindent=\parindent\ltextindent}
\font\Titel   = cmr10 scaled\magstep1
\centerline{\Titel References}

\bigskip \bigskip
{\parindent=1.3cm
\litem{[B/S]}B. Balcar and P. Simon, Cardinal invariants in Boolean
       spaces, General Topology and its Relations to Modern Analysis and
       Algebra V, Proc. Fifth Prague Topol. Symp. Symp. 1981, ed. J. Novak,
       Heldermann Verlag, Berlin, 1982, 39-47.
\litem{[B]}
       J.E. Baumgartner, Iterated forcing, Surveys of set theory
       (A.R.D. Mathias, editor), London Mathematical Society
       Lecture Note Series, no.87, Cambridge University Press,
       Cambridge, 1983, 1-59.
\litem{[B/Sp]}
       J.E. Baumgartner and O. Spinas, Independence
       and consistency proofs in quad\-rat\-ic form theory,
       Journal of Symbolic Logic, vol.57, no.4, 1991, 1195-1211.
\litem{[B/G]}
       W. Baur and H. Gross, Strange inner product spaces, Comment.
       Math. Helv. 52, 1977, 491-495.
\litem{[B/Sh]}
       A. Blass and S. Shelah, There may be simple $P_{\aleph _1}$-
       and $P_{\aleph _2}$-points and the Rudin-Keisler ordering may
       be downward directed, Annals of Pure and Applied Logic, vol.53,
       1987, 213-243.  
\litem{[vD]}
       E.K. van Douwen, The integers and topology, Handbook of
       set-theoretic topology, K. Kunen and J.E. Vaughan (editors),
       North-Holland, Amsterdam, 1984, 111-167.
\litem{[F]}
       D.H. Fremlin, Cicho\'n's diagram, Initiation \`a l'Analyse,
       Universit\'e Pire et Marie Curie, Paris, 1985.
\litem{[G]}
       H. Gross, Quadratic forms in infinite dimensional vector spaces,
       Progress in Mathematics, vol.1, Birkh\"auser, Boston, 1979.
\litem{[G/O]}
       H. Gross and E. Ogg, Quadratic spaces with few isometries,
       Comment. Math. Helv. 48, 1973, 511-519.
\litem{[Go]}
       M. Goldstern, Tools for your forcing construction, in: Set
       theory of the reals, Proceedings of the Bar Ilan Conference in
       honour of Abraham Fraenkel 1991, H. Judah (editor), 305-360.
\litem{[J]}
       T. Jech, Set theory, Academic Press, New York, 1978.
\litem{[J/Sh]}
       H. Judah and S. Shelah, Souslin forcing, Journal of Symbolic
       Logic, vol.53, no.4, 1988, 1188-1207. 
\litem{[K]}
       K. Kunen, Set theory. An introduction to independence proofs,
       North Holland, Amsterdam, 1980.
\litem{[Sh1]}
       S. Shelah, Proper forcing, Lecture Notes in Mathematics, vol.942,
       Springer.
\litem{[Sh2]}
       S. Shelah, On cardinal invariants of the continuum, Proceedings
       of the 6/83 Boulder conference in set theory, ed. J. Baumgartner,
       D. Martin and S. Shelah, Contemporary mathematics, vol.31, AMS,
       1984, 183-207.
\litem{[Sh3]}
       S. Shelah, Vive la difference I, Nonisomorphism of ultrapowers of
       countable models, in: Set theory of the continuum, ed. H. Judah, 
       W. Just, H. Woodin, Springer, New York, 1992, 357-405.
\litem{[Sh4]}
       S. Shelah, Strong negative partition relations below the
       continuum, Acta Math. Hungar. 58, no. 1-2, 1991, 95-100.
\litem{[Sh5]}
       S. Shelah, There are Jo\'nsson algebras in many inaccessible
       cardinals, in: Cardinal Arithmetic, Oxford University Press, 1994. 
\litem{[Sh6]}
       S. Shelah, Further cardinal arithmetic, Israel Journal of
       mathematics, in press ([Sh430] in Shelah's list of publications).
\litem{[Sh7]}
       S. Shelah, Colouring and $\aleph _2$-c.c. not productive, in 
       preparation ([Sh572] in Shelah's list of publications).
\litem{[Sp1]}
       O. Spinas, Konsistenz- und Unabh\"angigkeitsresultate in der
       Theorie der qua\-dra\-ti\-schen Formen, Dissertation, University of
       Z\"urich, 1989.
\litem{[Sp2]}
       O. Spinas, Iterated forcing in quadratic form theory, Israel
       Journal of Mathematics 79, 1991, 297-315.
\litem{[Sp3]}
       O. Spinas, An undecidability result in lattice theory,
       Abstracts of papers presented to the AMS, vol.11, no.2. p.161,
       March 1990.
\litem{[Sp4]}
       O. Spinas, Cardinal invariants and quadratic forms, in: Set
       theory of the reals, Proceedings of the Bar Ilan Conference in
       honour of Abraham Fraenkel 1991, H. Judah (editor), 563-581. 

}

\bigskip

{\bf Addresses:}  {\it First author:}  Department of Mathematics, The
Hebrew University of Jerusalem, Givat Ram, Jerusalem, ISRAEL.

\smallskip

{\it Second author:}  Department of Mathematics, University of California,
Irvine, CA 92717, USA

\bye